\documentstyle[twoside,12pt]{article}

\newcommand{\chapter}{\section}

\begin{document}

\newtheorem{Thm}{Theorem}
\newtheorem{Ax}{Axiom}
\newtheorem{Prop}{Proposition}
\newtheorem{Cor}[Prop]{Corollary}
\newtheorem{Main}{}
\renewcommand{\theMain}{}
\newtheorem{Lem}[Prop]{Lemma}
\newtheorem{Fact}{Fact}
\renewcommand{\theFact}{}

\newtheorem{Def}{Definition}
\newtheorem{rmk}{Remark}
\newenvironment{Rmk}{\begin{rmk}\em}{\end{rmk}}
\newtheorem{exm}{Example}
\newenvironment{Exm}{\begin{exm}\em}{\end{exm}}

\newcommand{\qed}{\par {\EM QED} }
\newtheorem{prf}{Proof}
\renewcommand{\theprf}{}
\newenvironment{Prf}{\begin{prf}\em}{\qed\end{prf}}
\newtheorem{prff}{}
\renewcommand{\theprff}{}
\newenvironment{Prff}{\begin{prff}\em}{\qed\end{prff}}





\newcommand{\YES}[1]{#1}
\newcommand{\NOT}[1]{}

\newcommand{\cA}{{\cal A}}
\newcommand{\cB}{{\cal B}}
\newcommand{\cC}{{\cal C}}
\newcommand{\cD}{{\cal D}}
\newcommand{\cE}{{\cal E}}
\newcommand{\cF}{{\cal F}}
\newcommand{\cG}{{\cal G}}
\newcommand{\cH}{{\cal H}}
\newcommand{\cI}{{\cal I}}
\newcommand{\cJ}{{\cal J}}
\newcommand{\cK}{{\cal K}}
\newcommand{\cL}{{\cal L}}
\newcommand{\cM}{{\cal M}}
\newcommand{\cN}{{\cal N}}
\newcommand{\cO}{{\cal O}}
\newcommand{\cP}{{\cal P}}
\newcommand{\cQ}{{\cal Q}}
\newcommand{\cR}{{\cal R}}
\newcommand{\cS}{{\cal S}}
\newcommand{\cT}{{\cal T}}
\newcommand{\cU}{{\cal U}}
\newcommand{\cV}{{\cal V}}
\newcommand{\cW}{{\cal W}}
\newcommand{\cX}{{\cal X}}
\newcommand{\cY}{{\cal Y}}
\newcommand{\cZ}{{\cal Z}}

\newcommand{\bbb}[1]{{\mbox{\bf #1}}}

\newcommand{\bN}{\bbb{N}}
\newcommand{\bZ}{\bbb{Z}}
\newcommand{\bR}{\bbb{R}}
\newcommand{\bC}{\bbb{C}}
\newcommand{\bQ}{\bbb{Q}}
\newcommand{\bT}{\bbb{T}}

\newcommand{\noind}[1]{{\setlength{\parindent}{0cm} #1}}
\newcommand{\parsk}{\par\medskip}

\newcommand{\varend}{\end{document}}

\newcommand{\LP}{\left(}  \newcommand{\RP}{\right)}
\newcommand{\LQ}{\left[}  \newcommand{\RQ}{\right]}
\newcommand{\LB}{\left\{} \newcommand{\RB}{\right\}}
\newcommand{\LA}{\left\langle} \newcommand{\RA}{\right\rangle}
\newcommand{\LS}{\left|}  \newcommand{\RS}{\right|}
\newcommand{\LN}{\left\|} \newcommand{\RN}{\right\|}
\newcommand{\bLP}{\bigl(}  \newcommand{\bRP}{\bigr)}

\newcommand{\ts}{{\otimes}} \newcommand{\TS}{{\bigotimes}}

\newcommand{\es}{\emptyset}
\newcommand{\sm}{\setminus} \newcommand{\st}{\subset}
\newcommand{\eps}{\varepsilon}

\newcommand{\beware}{~$\Psi\!\Psi\!\Psi$~}

\newcommand{\EM}{\em\bf}

\newcommand{\BE}{\begin{equation}}  \newcommand{\EE}{\end{equation}}
\newcommand{\BER}{$$\begin{array}}  \newcommand{\EER}{\end{array}$$}
\newcommand{\head}[1]{
             \par\bigskip\noind{{\large\bf#1}\nopagebreak\par}\medskip}
\newcommand{\chead}[1]{\begin{center}
              \par\bigskip{\large\bf#1}\nopagebreak\par\medskip
              \end{center}}

\newcommand{\fillitem}{\L{\embox{}}\vspace{-.6cm}}

\newcommand{\dfrac}[2]{\displaystyle{\frac{#1}{#2}}}
\newcommand{\toprove}[1]{{\buildrel\mbox{?}\over#1}}
\newcommand{\DEF}{{\buildrel\mbox{\scriptsize\,\,def\,\,}\over=}}
\newcommand{\BAR}[1]{\overline{#1}}
\newcommand{\I}{\infty}
\newcommand{\al}{\alpha} \newcommand{\be}{\beta}
\newcommand{\OM}{\Omega} \newcommand{\om}{\omega}
\newcommand{\la}{\lambda}
\newcommand{\mod}{\mbox{ mod }}
\newcommand{\half}{\frac{1}{2}}
\newcommand{\NI}{{n\to\infty}}
\newcommand{\ang}{{<\!\!\!)}}
\newcommand{\arc}[1]{\breve{#1}}

\newenvironment{txeqn}[1]{\begin{equation}\label{#1}\begin{array}{l}}{
                         \end{array}\end{equation}}



\newcommand{\id}{{\mbox{id}\,}}
\newcommand{\mns}{{(-1)^{\LS I\RS - \LS S\RS}}}

\title{Why Do Partitions Occur in Fa\`a di Bruno's Chain Rule
For Higher\NOT{Fr\'echet} Derivatives?}

\author{Eliahu Levy\\
Department of Mathematics\\
Technion -- Israel Institute of Technology,
Haifa 32000, Israel\\
email: eliahu@techunix.technion.ac.il}

\date{}

\maketitle

\begin{abstract}
It is well known that the coefficients in Fa\`a di Bruno's chain rule
for $n$-th derivatives of functions of one variable can be expressed
via counting of partitions\NOT{(see \cite{Johnson})}.
It turns out that this has a natural form as a formula\NOT{$(\ref{Main}.1)$}
for the vector case. Viewed as a purely algebraic fact, it is briefly
``explained'' in the first part of this note
why a proof for this formula leads to partitions. In the rest of the note
a proof for this formula is presented ``from first principles''
for the case of $n$-th Fr\'echet derivatives of mappings
between Banach spaces.
Again the proof ``explains'' why the formula has its form involving
partitions.\NOT{This formula relates to facts mentioned in \cite{B}.}
\end{abstract}


\setcounter{section}{-1}

\section{Introduction} \label{Intr}
If $I$ is a finite set, $\LS I\RS$  will denote its cardinality.
$\bN$ will denote the set of non-negative integers. $\cP(I)$ is
the power set of $I$.\parsk

Fa\`a di Bruno's formula (see \cite{Johnson} for an excellent
survey and bibiliography) gives a somewhat complicated expression
for $(g\circ f)^{(n)}(x)$, $f$ and $g$ functions of one variable,
in terms of derivatives up to order $n$ of $f$ at $x$ and of $g$
at $f(x)$. The coefficients can be expressed using numeration of
partitions of a set of $n$ elements (see below)\NOT{\cite}.
\parsk

It might seem that this settles the problem for mappings between
multi-dimensional spaces $X$: the existence of the $n$-th derivative
for the composite mapping follows from $1$-st derivative theorems,
the $n$-th derivatives in the vector case being treated as
iterated $1$-st derivatives, and the formula for it is obtained from
Fa\`a di Bruno's formula for one variable by composing with
linear mappings $X\to\bR$ or $\bR\to X$. Moreover, since it is clear, by
iterating the $1$-st derivative chain rule, that $(g\circ f)^{(n)}(x)$
is a certain polynomial in $f^{(k)}(x)$ and $g^{(k)}(f(x))$,
$1\le k\le n$, the problem is only to find its form, a purely algebraic
problem, and in this way it is treated in the literature.
\parsk

It turns out, however, that writing the formula for mappings between
vector spaces gives it a very natural form, obscured by the particulars
of the $1$-dimensional case.
This chain-rule formula for $n$-th derivatives of
mappings between vector spaces is constructed in terms of partitions
of a set $I$ of cardinality $n$, as follows:
\parsk

$$\LA(g\circ f)^{(n)}(x),\TS_{i\in I} v_i\RA =
\sum_{\pi\in\cH(I)}
\LA g^{(\LS\pi\RS)}(f(x)),\TS_{S\in\pi}\LA f^{(\LS S\RS)}(x),
\TS_{i\in S} v_i\RA\RA
\eqno(\ref{Main}.1)
$$
Where $f:X\to Y$ , $g:Y\to Z$ , $X$, $Y$, $Z$ are vector spaces,
$x, v_i\in X$,
$I$ is a set with $n$ elements and $\cH(I)$ is the set of partitions of
$I$, i.e.\ $\cH(I)$ is the subset of $\cP(\cP(I))$ consisting of all
disjoint collections of subsets of $I$ whose union is $I$ and
that do not contain the empty set.\NOT{$\LS\,\RS$ denotes number of
elements of a set.} Thus, for example:
$$\begin{array}{ll}
\cH(\es)=&\{\es\}\\
\cH(\{1\})=&\{\{1\}\}\\
\cH(\{1,2\})
=&\{\, \{\{1\},\{2\}\}\, , \,\{\{1,2\}\}\, \}\\
\cH(\{1,2,3\})=&\{\, \{\{1\},\{2\},\{3\}\}\, , \,\{\{1\},\{2,3\}\}\, ,
\,\{\{2\},\{1,3\}\}\, , \,\{\{3\},\{1,2\}\}\, ,
\{\{1,2,3\}\}\, \}
\NOT{\\  &\phantom{\{\,} \{\{1,2,3\}\}\, \}}
\end{array}$$
\parsk

For the one-dimensional case $X=Y=Z=$ the scalars
it suffices to let $v_i=1$ and $(3)$ reads:

$$\LA(g\circ f)^{(n)}(x), 1\RA =
\sum_{\pi\in\cH(I)}
\LA g^{(\LS\pi\RS)}(f(x)),\TS_{S\in\pi}\LA f^{(\LS S\RS)}(x),
1\RA\RA,
$$
that is:
$$(g\circ f)^{(n)}(x)=
\sum_{\pi\in\cH(I)}
g^{(\LS\pi\RS)}(f(x))\prod_{S\in\pi}f^{(\LS S\RS)}(x),
\eqno(\ref{Main}.1')
$$
a form of Fa\`a di Bruno's formula.
\parsk

If $(\ref{Main}.1)$ is viewed as a purely algebraic statement, one can deduce it in a
straightforward, though abstract, manner (see \S\ref{Alg}),
where the role of partitions is ``explained''.\parsk

In Theorem 1.\ in \S\ref{Main}, though, $(\ref{Main}.1)$ is proved
``from basic principles'' for Fr\'echet derivatives between Banach spaces.
Once guessed, $(\ref{Main}.1)$ can be proved in a straightforward manner using
induction. We have preferred to present a different
proof which explains directly the form of the formula. To this
end the derivatives are treated via $n$-th iterated differences,
i.e.\ alternating signs sums over vertices of $n$-dimensional
parallelepipeds, rather than iteratively. The formula is a consequence
of an identity (Lemma \ref{Lem2}) involving such sums. This identity
and the proof are best expressed in the language of free linear
spaces and free commutative algebras, which is therefore introduced.
Also, this approach to the $n$-th derivative allows one to require just
their existence in the strict sense (a definition is given in%
\NOT{the preliminaries section} \S\ref{Prel})
at the particular relevant points.
A slight complication arises from the fact that the identity of
Lemma \ref{Lem2} involves sums over parallelepipeds of dimensions
higher than $n$, the order of the highest participating derivative.
This is dealt with using Lemma \ref{Lem1}.
\parsk

In Bourbaki, Vari\'et\'es (\cite{Bourbaki})
the linear space of point distributions
of order $\le n$ at a point in a $C^{(n)}$-manifold is defined and its
properties stated.  The fact that this is properly defined and is
a functor may be proved using formula $(\ref{Main}.1)$.
\NOT{  (see remark after the proof of Theorem 1).}
Thus one may say that in spirit, $(\ref{Main}.1)$ is present in
\cite{Bourbaki}.
\parsk

\section[Where Do Partitions Come From? -- Algebraic Approach]%
{Where Do Partitions Come From? -- A Purely Algebraic Approach}
\label{Alg}
In this section the scalars $K$ may be any field. We insist on avoiding
division by integers, thus the scalar field may have any characteristic.
\parsk

Let us be abstract. The mappings will be between ``germs of $n$-manifolds
at a point'', in short $n$-germs, defined by their function algebra, namely
$\cA=\cA_n=K[[T_1,\ldots,T_n]]$ (formal power series), where the $T_j$ are
indeterminates. $\cA$ has a linear topology, taking as basic neighborhoods at $0$ the
positive integer powers of the ideal consisting of the formal power series
without constant term. Smooth mappings $f$ from the $n$-germ to the $m$-germ
are defined by continuous homomorphisms $\phi_f:\cA_m\to\cA_n$.
(The ``actual'' $m$ components of $f$ as a mapping are the ``functions''
$\phi_f(T_i),\;i=1,\ldots,m$.) The continuity of $\phi_f$ means just that any
coefficient of $\phi(F)$ is only a function of a finite number of
coefficients of $F$. Note that, by the continuity, since in the linear topology
of $\cA$\quad $T_i^k\to_{k\to\I}0$
also $\LP\phi_f(T_i)\RP^k\to_{k\to\I}0$, hence $\phi_f(T_i)$ have no
constant term -- indeed $f(0)=0$.
\parsk

Let $\cA'$ be the ``dual'', i.e.\ the $K$-vector space of continuous
functionals $\xi:\cA\to K$, $K$ given the discrete topology. Continuity
again means that $\xi(F)$\NOT{(to be denoted also by $\LA\xi,F\RA$)} depends
only on a finite set of coefficients of $F$.
\parsk

The multiplication on $\cA$ induces a {\EM comultiplication}
$\Delta:\cA'\to\cA'\TS\cA'$ by
$$(\Delta(\xi))(F\ts G):=\xi(FG)\qquad F,G\in\cA.\eqno(\ref{Alg}.1)$$
This comultiplication, with the counit $\xi\mapsto\xi(1)$, turns
$\cA'$ into a (coassociative and cocommutative) coalgebra over $K$.
$\cA'$ contains the Dirac delta $\delta$ at $0$, defined by $\delta(F):=$
the constant term in $F$, which satisfies $\Delta(\delta)=\delta\ts\delta$. 
We also have the primitive elements in $\cA'$, i.e.\ those satisfying
$\Delta(\xi)=\xi\ts\delta+\delta\ts\xi$, which means by $(\ref{Alg}.1)$
that the functional $\xi$ is a ``tangent vector'' at $0$. Indeed, they
constitute a vector space over $K$ isomorphic to $K^n$ by the isomorphism
$\xi\mapsto(\xi(T_1),\ldots,\xi(T_n))$. Denote this ``tangent space'' of the
$n$-germ by $V=V_n$. We identify a ``vector'' in $V$ with the corresponding
``vector'' in $K^n$.
\parsk

But there is more structure: one may define the {\EM convolution} of two
elements $\xi,\eta$ of $\cA'$ (giving an element of $\cA'$), and the
convolution of an element $\xi$ of $\cA'$ and a $F\in\cA$
(giving an element of $\cA$) by:
$$(\xi*\eta)(F):=\xi(T\mapsto\eta(S\mapsto F(S+T)))\qquad
\xi*F:=(T\mapsto\xi(S\mapsto F(S+T)))\qquad F\in\cA.\eqno(\ref{Alg}.2)$$
The first convolution makes $\cA'$ also into an (associative and commutative)
algebra, with unit $\delta$, which acts linearly on $\cA$ by the second
convolution. In particular, as one easily sees, the primitive elements
$v\in V$ act on $\cA$ as the directional partial derivatives with direction
$v$. Indeed, one may write:
$$\LA F',v\RA=v*F\qquad v(F)=\LA F'(0),v\RA\qquad v\in V\,,\,F\in\cA.$$ 
Convolutions of elements $v$ in $V$ will act, by convolution with members
of $\cA$, as composition of the actions of the $v$, i.e.\ by higher-order
differential operators. One finds that one may also write:
$$(v_1*\cdots*v_k)(F)=\LA F^{(k)}(0),v_1\ts\cdots\ts v_k\RA\qquad
v_i\in V\,,\,F\in\cA.\eqno(\ref{Alg}.3)$$
Indeed:
\begin{eqnarray*}
\LA F^{(k)}(0),v_1\ts\cdots\ts v_k\RA=
\LA F^{(k)},v_1\ts\cdots\ts v_k\RA|_0=\\
=(v_1*\cdot*v_n*F)|_0=(v_1*\cdots*v_n*\delta)(F)=
(v_1*\cdots*v_n)(F)
\end{eqnarray*}
\parsk

Yet there is a difference between the comultiplication and the convolution
multiplication: since any smooth mapping $f$ defines a homomorphism $\phi_f$
of the algebra $\cA$, it will induce a dual map on $\cA'$ (which we again
denote by $\phi_f$) which will preserve the comultiplication, i.e.\
$\phi_f(\Delta(\xi))=\Delta(\phi_f(\xi))$, (Applying
a mapping such as $\phi_f$ to the tensor product is understood factorwise.)
Thus, the coalgebra structure is an invariant of the manifold structure,
so to speak. On the other hand the convolution was defined in $(\ref{Alg}.2)$ using
the additive structure of $K^n$, so to speak, and a smooth mapping
need preserve it only if it is linear. Nevertheless, using the fact that the
convolution in $\cA'$ is, in a sense, the map $\cA'\to\cA'\ts\cA'$ induced by
the smooth mapping $K^n\times K^n\to K^n$ given by addition, one shows that
convolution preserves the comultiplication, hence $\cA'$ is a bialgebra,
(it turns out to be a Hopf algebra).
\parsk

This can be used to compute $\Delta$ on $*_{i\in I}v_i$,\,
$v_i\in V$,\, $I$ a finite set.
(Recall that $\delta$ is the unit of convolution):
$$\Delta\LP*_{i\in I}v_i\RP=*_{i\in I}\Delta(v_i)=
*_{i\in I}(v_i\ts\delta+\delta\ts v_i)=
\sum_{S\st I}\LP*_{i\in S}v_i\RP\ts\LP*_{i\in I\sm S}v_i\RP.
\eqno(\ref{Alg}.4)$$
\parsk

One defines $\Delta^{(3)}:\cA'\to\cA'\ts\cA'\ts\cA'$ by
$\Delta^{(3)}:=(\Delta\ts\id)\circ\Delta=(\id\ts\Delta)\circ\Delta$
(the latter being equal by coassociativity), and one has 
$\LQ\Delta^{(3)}(\xi)\RQ(F\ts G\ts H)=\xi(FGH)$. Similarly
for $\Delta^{(n)}$,\,$n\in\bN$. For $v\in V$ one has
$\Delta^{(3)}v=v\ts\delta\ts\delta+\delta\ts v\ts\delta+
\delta\ts\delta\ts v$ etc. In analogy to $(\ref{Alg}.4)$, one has,
for $v_i\in V$,\, $I$ a finite set:
$$\Delta^{(k)}\LP*_{i\in I}v_i\RP=\sum_{I=\mbox{disj.}\,\cup_{j=1}^kS_j}
\TS_{j=1}^k*_{i\in S_j}v_i.\eqno(\ref{Alg}.5)$$ 
\parsk

Let $f$ be a smooth mapping from the $n$-germ to the $m$-germ.
Applying $(\ref{Alg}.3)$ for $F=f_j=\phi_f(T_j)$ one obtains
(recall that $f^{(k)}(0)$ is a linear map $S^k(V_n)\to V_m$
-- for the notation $S^k(V)$ see \S\ref{SS:ST}):
$$\LQ\phi_f(v_1*\cdots*v_k)\RQ(T_j)=
\LQ\LA f^{(k)}(0),v_1\ts\cdots\ts v_k\RA\RQ(T_j)=
\LQ\LA f^{(k)}(0),v_1\ts\cdots\ts v_k\RA\RQ_j
\quad v_i\in V\quad j=1,\ldots,m.
\eqno(\ref{Alg}.6).$$
We claim that, for $v_i\in V$,\, $I$ a finite set (for the notation $\cH(I)$
see \S\ref{Intr}\NOT{or \S\ref{Part}}):
$$\phi_f(*_{i\in I}v_i)=\sum_{\pi\in\cH(I)}
*_{S\in\pi}\LA f^{(\LS S\RS)}(0),\ts_{i\in S}v_i\RA.\eqno(\ref{Alg}.7)$$
To get $(\ref{Alg}.7)$, it suffices that both sides give the same value when
applied to a monomial on the $T_j$'s, that is, that applying $\Delta^{(k)}$
to both sides one obtains elements of the $k$-th tensor power of $\cA'$ that
give the same value when applied to tensors of $T_j$'s. But that follows in a
straightforward (though somewhat clumsy) manner from $(\ref{Alg}.5)$ and $(\ref{Alg}.6)$ and
the fact that $\phi_f$ preserves the comultiplication.\parsk

Now $(\ref{Main}.1)$ follows from $(\ref{Alg}.7)$ and $(\ref{Alg}.6)$:
indeed, if $f$ and $g$ are smooth mappings, $v_i\in V$,\, $I$ a finite set
and $j=1,\ldots,n$ where $n$ is the dimension of the image germ of $g$:
\begin{eqnarray*}
\LQ\LA(g\circ f)^{(\LS I\RS)}(0),\TS_{i\in I} v_i\RA\RQ_j =
\LQ\phi_g(\phi_f(*_{i\in I}v_i))\RQ(T_j)=\\
=\LQ\phi_g\LP\sum_{\pi\in\cH(I)}
*_{S\in\pi}\LA f^{(\LS S\RS)}(0),\TS_{i\in S}v_i\RA\RP\RQ(T_j)=\\
=\LQ\sum_{\pi\in\cH(I)}
\LA g^{(\LS\pi\RS)}(0),\TS_{S\in\pi}\LA f^{(\LS S\RS)}(0),
\TS_{i\in S} v_i\RA\RA\RQ_j
\end{eqnarray*}

\section{Preliminaries to the Banach Space Approach} \label{Prel}
Let $X$, $Y$, $Z$ denote linear spaces over a field (the scalar field)
-- from \S\ref{Main} on Banach spaces unless otherwise stated.
The choice of $\bR$ or $\bC$ as scalar field is insignificant and will
not be mentioned.
\parsk

\subsection{Symmetric Tensors} \label{SS:ST}
We consider the {\EM symmetric power} $S^n(X)$ of $X$, defined as the
quotient space of the n-fold (algebraic) tensor product $\TS^n X$ with
respect to its subspace spanned by the elements
$$x_1\ts x_2\ts\dots\ts x_n
- x_{\sigma(1)}\ts x_{\sigma(2)}\ts\dots\ts x_{\sigma(n)}$$
where $x_i\in X$ and $\sigma$ is a permutation of $\{1,2,\dots,n\}$.
The image of a tensor  $x_1\ts\dots\ts x_n$ in $S^n(X)$ by the
quotient map will still be
denoted by  $x_1\ts\dots\ts x_n$.  Thus, $(x_1,\dots,x_n)\mapsto
x_1\ts\dots\ts x_n \in S^n(X)$
is multilinear and symmetric and one may freely employ the notation
$\TS_{i\in I} x_i$ for a finite family $(x_i)_{i\in I}$.
\parsk

In particular, $S^0(X)$ may be identified with the scalar field.
\parsk

In the direct sum $S(X):=\bigoplus_0^\I S^n(X)$,\quad$\ts$ is an associative
and commutative multiplication. $S(X)$ is called the {\EM symmetric algebra}
of $X$.
\parsk

When $X$ and $Y$ are Banach spaces, the bounded linear operators
$T: S^n(X)\to Y$, i.e.\ those satisfying
$$\LN T\RN :=
\sup_{\LN x_i\RN\le 1}\LS\LA T\,,\,x_1\ts\dots\ts x_n\RA\RS < \infty$$
form, with the above norm, a Banach space, to be denoted by
$\cL^{(n)}(X,Y)$.
Of course, this space may be identified with the Banach space of
bounded symmetric $n$-multilinear operators from $X^n$ to $Y$.
\parsk

\subsection {The $n$-th strict (Fr\'echet) derivative} It may be defined,
instead of iteratively,
via $n$-th order differences (or, in other words, alternating signs
sums over vertices of parallelepipeds), as follows:

\begin{Def} \label{Def1}
Let $X$, $Y$ be Banach spaces, $U\subset X$, $f:U\to Y$,
$x\in U^o$ (the interior of U).
Let $n\in\bN$.  We say that $f$  has a (strict) $n$-th (Fr\'echet)
derivative at $x$, the derivative being the element $f^{(n)}$ of
$\cL^{(n)}(X,Y)$, if $f^{(n)}$ satisfies, for a set $I$ with
$\LS I\RS =n$:
$$ \sum_{S\subset I}\mns f\LP\bar v+\sum_{i\in S} v_i\RP =
\LA f^{(n)}(x)\,,\,\TS_{i\in I}v_i\RA+
o_{\bar v\to x,v_i\to 0}\LP\prod_{i\in I}\LN v_i\RN\RP
\eqno(\ref{Prel}.1)$$
\end{Def}

Such an $f^{(n)}$ is necessarily unique, if it exists.
\parsk

Defined this way, $f^{(n)}(x)$  may exist even if previous derivatives
fail to exist in a neighborhood of $x$ or even at $x$ itself  (for
example any non-continuous linear operator $f$ has $f^{(n)}\equiv 0$ ,
$n\ge 2$).
\parsk

Note that $f^{(0)}(x) = y$  means that $f$ is continuous at $x$  and
$f(x) = y$.
\parsk

The proof of the following fact is standard:
\par\medskip
{\EM Fact:} Suppose $f^{(n)}(x)$ exists for every $x$ in an open
$U\subset X$. Then:
\begin{itemize}
\item[(i)] $f^{(n)}:U\to\cL^{(n)}(X,Y)$ is continuous
\item[(ii)] If $\bar v\in U$ , $v_i\in X$ are such that the
"parallelepiped"
$$\left\{\bar v+\sum t_i v_i : t=(t_1,\dots,t_n)\in [0,1]^n\right\}
\subset U$$
then ($\LS I \RS = n$):
$$\sum_{S\subset I} \mns f\LP\bar v+\sum_{i\in I} v_i \RP =
\int_{t\in [0,1]^n} \LA f^{(n)}\LP\bar v + \sum_{i\in I}t_i v_i\RP\,,
\,\TS_{i\in I} v_i\RA\,dt
$$
$dt$ being the n-dimensional Lebesgue measure.\parsk
\item[(iii)] (consequence of (ii) ):  If $x\in U$ , $m\in\bN$ then
$f^{(n)}$ has a (strict) $m$-th derivative at $x$ iff
$f$ has a (strict) $m+n$ -th derivative at $x$
and then:
$$\LA f^{(n+m)}(x),\TS_{1\le i\le m+n} v_i\RA =
\LA\LA\LP f^{(n)}\RP^{(m)}(x),
\TS_{1\le i\le m}v_i\RA,\TS_{m+1\le i\le m+n}v_i\RA.$$
\end{itemize}

(iii) implies that the above definition of $f^{(n)}(x)$ coincides with the
iterative definition whenever $f^{(m)}$ , $m<n$ exist in a neighborhood
of $x$.
\parsk

\section{The Chain Rule for the $n$-th Fr\'echet Derivative} \label{Main}
\begin{Main} {\EM Theorem 1.}  \label{Tm:Main}
Let $n\in\bN$ , $f:U\subset X\to Y$ ,
$g:V\subset Y\to Z$ ($X$,$Y$,$Z$ Banach spaces) be
such that $f(U)\subset V$. Let $x\in U^o$ (interior of U), be such
that $f^{(k)}$ , $0\le k\le n$
exist at $x$, $f(x)\in V$ and $g^{(k)}$ , $0\le k\le n$ exist at
$f(x)$. Then 
$(g\circ f)^{(n)}(x)$ exists and we have the formula ($\LS I \RS=n$):
$$\LA(g\circ f)^{(n)}(x),\TS_{i\in I} v_i\RA =
\sum_{\pi\in\cH(I)}
\LA g^{(\LS\pi\RS)}(f(x)),\TS_{S\in\pi}\LA f^{(\LS S\RS)}(x),
\TS_{i\in S} v_i\RA\RA
\eqno(\ref{Main}.1)
$$
All derivatives are in the sense of Definition \ref{Def1}.
\end{Main}

{\EM Proof.} Let us remark first, that if it is assumed that the
derivatives of order $<n$ exist in neighborhoods of $x$ and
$f(x)$ and one employs the iterative definition, then $(\ref{Main}.1)$ may
be proved by induction in a straightforward manner (assuming
one guesses the formula in advance). We shall rather give
a proof that assumes existence of the derivatives only at
$x$ and $f(x)$ using Definition \ref{Def1}, and moreover "explains"
why $(\ref{Main}.1)$ has its form.
\parsk

Before we proceed with the proof, let us introduce the following
auxiliary concepts and notations:
\parsk

Assume, for the moment, that $X$ is just a set. Denote by $\cE X$
the free linear space with basis denoted by $\LP\cE x\RP_{x\in X}$.
In this way, $\cE$ is a functor: for any sets $X$, $Y$ and function
$f:X\to Y$ we have a linear function $f^\cE: \cE X\to \cE Y$
defined by $f^\cE\LP\cE x\RP = \cE\LP f(x)\RP$. In case $Y$ is a linear
space, one has also a linear $f^{L}: \cE X\to Y$ defined by
$f^{L}\LP\cE x\RP = f(x)$.
\parsk

For a linear space $X$, $\cE X$ is a commutative algebra where
$\cE x\cdot\cE y = \cE(x+y)$ , i.e. $X$ is the group-algebra of
the additive
group $X$. It has the unit element $1=\cE 0$. For linear $f:X\to Y$,
$f^\cE$ is an algebra homomorphism preserving unit element.
\parsk

$\cE$ may be iterated: for any set $X$ we have the commutative
algebra $\cE\cE X$.
\parsk

For the case that $X$,$Y$,$Z$ are linear spaces and $f:X\to Y$ ,
$g:Y\to Z$ any functions, we have the following formulas:
   
$$ f(x) = f^L\LP\cE x\RP \eqno(\ref{Main}.2)$$
$$ \cE f(x) = \cE f^L\LP\cE x\RP = f^{L\cE}\LP\cE\cE x\RP \eqno(\ref{Main}.2')$$
$$ gf(x) = g^L\LP\cE f(x)\RP = g^Lf^{L\cE}\LP\cE\cE x\RP \eqno(\ref{Main}.2'')$$
Note that $f^{L\cE}$ is always an algebra homomorphism.
\parsk

Let us rewrite the left-hand side expression of $(\ref{Prel}.1)$
(definition \ref{Def1}) using these notations:
$$\begin{array}{l}
\sum_{S\subset I}\mns f\LP\bar v +\sum_{i\in S}v_i\RP
= \sum_{S\subset I}\mns f^L\LP\cE\LP\bar v +\sum_{i\in S}v_i\RP\RP\\
=f^L\LP\sum_{S\subset I}\mns\cE\bar v\prod_{i\in S}\cE v_i\RP
=f^L\LP\cE\bar v\prod_{i\in I}\LP\cE v_i - 1\RP\RP.
\end{array}$$
Thus $(\ref{Prel}.1)$ takes the form:
$$f^L\LP\cE\bar v\prod_{i\in I}\LP\cE v_i - 1\RP\RP
=\LA f^{(n)}(x),\TS_{i\in I}v_i\RA +
o_{\bar v\to x, v_i\to 0}\LP\prod_{i\in I}\LN v_i\RN\RP
\eqno(\ref{Main}.3)
$$
\begin{Lem} \label{Lem1}
If $f^{(n)}$ exists but one takes $\LS I\RS>n$,
then for any $J\subset I$, $\LS J\RS=n$, the left hand-side of $(\ref{Main}.3)$ is
$o_{\bar v\to x, v_i\to 0}\LP\prod_{i\in J}\LN v_i\RN\RP$.
\end{Lem}

\begin{Prff} {\EM  Proof of Lemma \ref{Lem1}.}
$$\begin{array}{l}
f^L\LP\cE\bar v \prod_{i\in I}\LP\cE v_i-1\RP\RP
=f^L\LP\cE\bar v \prod_{i\in I\sm J}\LP\cE v_i-1\RP
\prod_{i\in J}\LP\cE v_i-1\RP\RP= \\
=f^L\LP\cE\bar v \sum_{S\subset I\sm J}
(-1)^{\LS I\sm J\RS-\LS S\RS} \prod_{i\in S}\cE v_i
\prod_{i\in J}\LP\cE v_i-1\RP\RP= \\
=\sum_{S\subset I\sm J} (-1)^{\LS I\sm J\RS-\LS S\RS}
f^L\LP\cE\LP\bar v+\sum_{i\in S}v_i\RP \prod_{i\in J}\LP\cE v_i-1\RP\RP
\end{array}$$
By $(\ref{Main}.3)$ (note that $\LS J\RS=n$), the last expression is
$$
\LP\sum_{S\subset I\sm J} (-1)^{\LS I\sm J\RS-\LS S\RS}\RP
\LA f^{(n)}(x),\TS_{i\in J}v_i\RA +
o_{\bar v\to x, v_i\to 0}\LP\prod_{i\in J}\LN v_i\RN\RP
$$
And the lemma follows from the fact that the first term vanishes, since
$$\sum_{S\subset I\sm J} (-1)^{\LS I\sm J\RS-\LS S\RS}=
(1-1)^{\LS I\sm J\RS}=0$$ because $\LS I\sm J\RS>0$.\parsk
\end{Prff}

Now $(\ref{Main}.1)$ will follow, using $(\ref{Main}.3)$, from the following
identity which will be proved later:

\begin{Lem} \label{Lem2}
For any finite set $I$ and vectors $\bar v, v_i$ in a linear space $X$: 
$$
\sum_{S\st I}\mns\cE\cE\LP\bar v+\sum_{i\in S}v_i\RP =
\sum_{\al\st\cP(I),\:\cup\al=I,\:\es\notin\al}
\cE\cE\bar v\prod_{A\in\al}\LP\cE\LQ\cE\bar v
\prod_{i\in A}(\cE v_i-1)\RQ-1\RP 
\eqno(*)
$$
\end{Lem}
\par\medskip

Proceeding with the proof of $(\ref{Main}.1)$, one has, by $(\ref{Main}.2)$, $(\ref{Main}.2')$ and $(\ref{Main}.2'')$,
using $(*)$:

$$\begin{array}{l}
\sum_{S\st I}\mns g\LP f\LP\bar v+\sum_{i\in I}v_i\RP\RP = \\   
= g^Lf^{L\cE}\LP\sum_{S\st I}\mns\cE\cE\LP\bar v+\sum_{i\in S}v_i\RP\RP = \\
= g^Lf^{L\cE}\LP\sum_{\al\st\cP(I),\:\cup\al=I,\:\es\notin\al}
\cE\cE\bar v\prod_{A\in\al}\LP\cE\LQ\cE\bar v
\prod_{i\in A}(\cE v_i-1)\RQ-1\RP\RP = \\
=\sum_{\al\st\cP(I),\:\cup\al=I,\:\es\notin\al}
g^L\LP\cE f(\bar v)\prod_{A\in\al}\LP\cE f^L\LQ\cE\bar v 
\prod_{i\in A}(\cE v_i-1)\RQ-1\RP\RP = 
\sum_{\al\st\cP(I),\:\cup\al=I,\:\es\notin\al}D_\al.
\end{array}$$

$D_\al$ is an expression of the form appearing in the left-hand side
of $(\ref{Main}.3)$, with $f$ replaced by $g$, $I$ by $\al$ and the vectors $\bar v$ and
$v_i$ by $\bar w=f(\bar v)$ and
$w_A=f^L\LP(\cE\bar v\prod_{i\in A}(\cE v_i-1)\RP$.

Since $f$ is continuous at $x$ (having a strict $0$-derivative) ,
$\bar w\to f(x)$ as $\bar v\to x$.  Also, for $\es\ne A\st I$,

$$\begin{array}{l}
w_A=f^L\LP(\cE\bar v\prod_{i\in A}(\cE v_i-1)\RP =
\LA f^{\LS A\RS}(x),\TS_{i\in A}v_i\RA + o\LP\prod_{i\in A}\LN v_i\RN\RP = 
O\LP\prod_{i\in A}\LN v_i\RN\RP
\end{array}$$
which tends to $0$ as $\bar v\to x,\: v_i\to 0$.

Now, since $g$ is assumed having strict derivatives at $f(x)$  up
to order  n , $(\ref{Main}.3)$ and Lemma \ref{Lem1} applied to $g$ give:

$$\begin{array}{l}
D_\al = g^L\LP\cE\bar w\prod_{A\in\al}(\cE w_A-1)\RP = \\
= \left\{\begin{array}{ll}
\LA g^{(\LS\al\RS)}(f(x)),\TS_{A\in\al}w_A\RA +
o\LP\prod_{A\in\al}\LN w_A\RN\RP =
O\LP\prod_{A\in\al}\LN w_A\RN\RP &\mbox{ if }\LS\al\RS\le n\\
o\LP\prod_{A\in J}\LN w_A\RN\RP
\mbox{ for any } J\st\al \mbox{ with } \LS J\RS=n, &\mbox{ if }\LS\al\RS>n
\end{array}\right.
\end{array}
\eqno(\ref{Main}.4)
$$
(All $O$'s and $o$'s are for $\bar v\to x,\: v_i\to 0$).

In order to deal with the second case in $(\ref{Main}.4)$, note that we have:

\begin{Lem} \label{Lem3}
If $\cup\al=I$, then $\exists\be\st\al$ such that $\cup\be=I$ and
$\LS\be\RS\le n$.
\end{Lem}

\begin{Prf}
Just pick for each $i\in I$ an $A\in\al$ such that $i\in A$ and let
$\be$ be the set of $A$'s picked.
\end{Prf}

Thus, if $\cup\al=I$ and $\LS\al\RS>n$ we can always find a $J\st\al$ with
$\LS J\RS=n$, $\cup J=I$ and
$D_\al = o\LP\prod_{A\in J}\LN w_A\RN\RP = o\LP\prod_{i\in I}\LN v_i\RN\RP$.
If $\cup\al=I$  and $\LS\al\RS\le n$ but $\al$ is not a partition, then
$D_\al=O\LP\prod_{A\in\al}\LN w_A\RN\RP = o\LP\prod_{i\in I}\LN v_i\RN\RP,$
since the $A$'s cover $I$ with redundancy. So we are left with the
$D_\al$'s for $\al\in\cH(I)$, for which we have:

$$\begin{array}{l}
D_\al=\LA g^{(\LS\al\RS)}(f(x)),\TS_{A\in\al} w_A\RA +
o\LP\prod_{w\in A}\LN w_A\RN\RP = \\
= \LA g^{(\LS\al\RS)}(f(x)) , \TS_{A\in\al}
\LQ\LA f^{(\LS A\RS)}(x),\TS_{i\in A}v_i\RA 
+ o\LP\prod_{i\in A}\LN v_i\RN\RP\RQ\RA
+ o\LP\prod_{i\in I}\LN v_i\RN\RP = \\
= \LA g^{(\LS\al\RS)}(f(x)),
\TS_{A\in\al}\LA f^{(\LS A\RS)}(x),\TS_{i\in A}v_i\RA\RA +
o\LP\prod_{i\in I}\LN v_i\RN\RP
\end{array}$$
which implies $(\ref{Main}.1)$.
\par\medskip

The only thing left to be done is the:

\par\medskip

\begin{Prff}{\EM Proof of Lemma \ref{Lem2}} (i.e. of $(*)$):

$$\begin{array}{l}
\sum_{S\st I}\mns\cE\cE\LP\bar v+\sum_{i\in S}v_i\RP=\\
=\sum_{S\st I}\mns\cE\LP\cE\bar v\prod_{i\in S}\cE v_i\RP=\\
=\sum_{S\st I}\mns\cE\LQ\cE\bar v\prod_{i\in S}\LP1+(\cE v_i-1)\RP\RQ=\\
=\sum_{S\st I}\mns\cE\LQ\cE\bar v\sum_{A\st S}\prod_{i\in A}(\cE v_i-1)\RQ=\\
=\sum_{S\st I}\mns\cE
\LQ\cE\bar v+\sum_{A\st S,\:A\ne\es}\cE\bar v\prod_{i\in A}(\cE v_i-1)\RQ=\\
=\sum_{S\st I}\mns\cE
\cE\bar v\prod_{A\st S,\:A\ne\es}\cE\LP\cE\bar v\prod_{i\in A}(\cE v_i-1)\RP
=\\
=\sum_{S\st I}\mns\cE\cE\bar v\prod_{A\st S,\:A\ne\es}
\LQ1+\LP\cE\LQ\cE\bar v\prod_{i\in A}(\cE v_i-1)\RQ-1\RP\RQ=\\
=\sum_{S\st I}\mns\cE\cE\bar v\sum_{\al\st\cP(S),\:\es\notin\al}
\prod_{A\in\al}\LP\cE\LQ\cE\bar v\prod_{i\in A}(\cE v_i-1)\RQ-1\RP=\\
=\sum_{\al\st\cP(I),\:\es\notin\al}
\LP\sum_{S\st I,\:\cup\al\st S}\mns\RP\cE\cE\bar v
\prod_{A\in\al}\LP\cE\LQ\cE\bar v\prod_{i\in A}(\cE v_i-1)\RQ-1\RP
\end{array}$$

and $(*)$ follows from the fact that for fixed $\al$
$$\sum_{S\st I,\:\cup\al\st S}\mns=(1-1)^{\LS I\RS-\LS\cup\al\RS}=
\left\{\begin{array}{ll}1 &\mbox{if }\cup\al=I\\0 &\mbox{otherwise}
\end{array}\right. .$$

\end{Prff}
\par\medskip
This completes the proof of Theorem 1.

\end{document}